\documentclass{amsart}

\usepackage{amsfonts}
\usepackage{amsmath}
\usepackage{amsthm}
\usepackage{amssymb}
\usepackage{textcomp}
\usepackage{cite}

\newtheorem{theorem}{Theorem}[section]

\newtheorem{corollary}[theorem]{Corollary}

\theoremstyle{definition}

\newtheorem{remark}[theorem]{Remark}

\newtheorem *{Theorem A}{Theorem A}
\newtheorem *{Theorem B}{Theorem B}

\begin{document}
\title[Graphs and graded algebras]{A graph-theoretic approach for comparing dimensions of components in
simply-graded algebras }
\author{Yuval Ginosar}

\address{Department of Mathematics, University of Haifa, Haifa 31905, Israel}
\email{ginosar@math.haifa.ac.il}
\author{Ofir Schnabel}
\email{os2519@yahoo.com}

\bibliographystyle{abbrv}

\begin{abstract}
Any simple group-grading of a finite dimensional complex algebra
induces a natural family of digraphs. We prove that $|E\circ
E^{\text{op}}\cup E^{\text{op}}\circ E|\geq |E|$ for any digraph
$\Gamma =(V,E)$  without parallel edges, and deduce that for any
simple group-grading, the dimension of the trivial component is
maximal.
\end{abstract}
\date{\today\vspace{-0.5cm}}
\maketitle

\section{Introduction}
Let $\Gamma =(V,E)$ be a finite digraph, where $V=\{v_1,v_2,\ldots ,v_n\}$.
Set $$D^-(v_i)=\{v_j\in V|(v_j,v_i)\in E\},\quad D^+(v_i)=\{v_j\in V|(v_i,v_j)\in E\}.$$
Two vertices $v_i,v_j\in V$ are \textit{mutually neighbored}
if $$D^-(v_i)\cap D^-(v_j)\neq \emptyset \quad  \text{or} \quad D^+(v_i)\cap D^+(v_j)\neq \emptyset.$$
Let $$T(\Gamma)=E\circ E^{\text{op}}\cup  E^{\text{op}}\circ E=\{(v_i,v_j)|v_i,v_j \text{ are mutually neighbored}\}.$$
For simplicity, we sometimes abuse notation and consider
$$T(\Gamma)=\{(i,j)|v_i,v_j \text{ are mutually neighbored}\}.$$
\begin{Theorem A}
For any finite digraph without parallel edges $\Gamma =(V,E)$,
$$|T(\Gamma)|\geq|E|.$$
\end{Theorem A}
A consequence of Theorem A is an estimation of the dimension of the homogeneous components in
simple group-graded algebras.
A \textit{grading} of an algebra $\Lambda$ by a group $G$ is a vector space decomposition
\begin{equation}\label{eq:algebragrading}
\Lambda=\bigoplus _{g\in G} \Lambda_g
\end{equation}
such that $\Lambda_g\Lambda_h\subseteq \Lambda_{gh}$.
An algebra $\Lambda$ is \textit{$G$-simple} with respect to a group grading~\eqref{eq:algebragrading}
if it admits no non-trivial graded ideals.
We prove the following theorem, which is implicit in \cite{giambruno}.
\begin{Theorem B}
Let $\Lambda$ be a complex $G$-simple algebra with respect to the grading~\eqref{eq:algebragrading}.
Denote the trivial element of $G$ by $e$.
Then dim$_{\mathbb{C}}(\Lambda_e)\geq$dim$_{\mathbb{C}}(\Lambda_g)$ for any $g\in G$.
\end{Theorem B}
Theorem A is proven in \S\ref{graph results}.
In \S\ref{gradings} we associate a natural digraph to any simple grading of a complex algebra
using a classification theorem due to Bahturin, Sehgal and Zaicev.
Then, using Theorem A we prove Theorem B.

\section{Proof of Theorem A}\label{graph results}
Let $V=\{v_1,v_2,\ldots ,v_n\}$. Obviously, we may assume that
$\Gamma$ has no isolated vertices. In particular, we may assume
that $(i,i)\in T(\Gamma)$ for any $1\leq i\leq n$. Another
assumption we can adopt is that $\Gamma$ admits no loops. Indeed,
for every loop $(v_i,v_i)\in E$, the condition that either
$(v_i,v_j)\in E$ or $(v_j,v_i)\in E$ says that $v_i\in
D^-(v_i)\cap D^-(v_j)$ or $v_i\in D^+(v_i)\cap D^+(v_j)$
respectively, which is the same as saying that both
$(i,j),(j,i)\in T(\Gamma)$. Hence the number of ordered pairs that
$v_i$ contributes to $T(\Gamma)$ is at least as the number of
edges that this vertex contributes to $E$.\\
The theorem is clear for $|V|=1,2$. We proceed by induction on the
number of vertices in $\Gamma $. Assume that the theorem holds for
digraphs with $|V|\leq n-1$. Let $\Gamma=(V,E)$ where $|V|=n$.
We distinguish between the Eulerian and the non-Eulerian cases:\\
\underline{\textbf{$\Gamma$ is Eulerian.}}\\
For every $1\leq i\leq n$ denote by $r_i=|D^-(v_i)|=|D^+(v_i)|$.
Let $k$ be the length of a shortest directed cycle in $\Gamma $.
Without loss of generality we may assume that the vertices in the cycle are $\{v_i\}_{1\leq i\leq k}$,
such that $v_i \in D^-(v_{i+1})$ for $1\leq i< k$ and $v_k\in D^-(v_1)$.
By minimality of the length of the cycle, if $i\not \equiv j+1($mod $k)$, then $v_j\not \in D^-(v_i)$.
Again we distinguish between two different cases.\\
\textbf{Case 1:} Assume that for any $1\leq i\neq j\leq k$ the ordered pair $(i,j)\not \in T(\Gamma)$.
By removing the vertices $v_1,v_2,\ldots ,v_k$ as well as their corresponding edges
we get a new graph $\Gamma ^{\shortmid}=(V^{\shortmid},E^{\shortmid})$ with $n-k$ vertices, which satisfies the
theorem by the induction assumption. That is,
\begin{equation}\label{eq:ind}
|T(\Gamma ^{\shortmid})|\geq |E^{\shortmid}|.
\end{equation}
The number of edges that were removed is
\begin{equation}\label{eq:edges}
|E|-|E^{\shortmid}|=\sum _{i=1}^k 2r_i-k.
\end{equation}
We count the number of ordered pairs that were removed from $T(\Gamma)$.
For $1\leq i<k$, let
$$ C_i=\{(v_i,v)|v\in D^-(v_{i+1})\}\cup \{(v,v_i)|v\in D^-(v_{i+1})\}(\subseteq T(\Gamma)),$$
and also
$$ C_k=\{(v_k,v)|v\in D^-(v_1)\}\cup \{(v,v_k)|v\in D^-(v_1)\}(\subseteq T(\Gamma)).$$
In order to compute the cardinality of $C_i$, for $1\leq i<k$, every
$v\in D^-(v_{i+1})$ is counted twice, except $v_i$ itself which is counted
only once. This argument, as well as a similar argument for $C_k$ yields
\begin{equation}\label{eq:C_i}
|C_i|= 2r_{i+1}-1 , 1\leq i< k ,\quad |C_k|= 2r_1-1.
\end{equation}
 The condition that for $1\leq i\neq l\leq k$ the
ordered pair $(i,l)\not \in T(\Gamma)$ ensures that for $i\neq l$,
$C_i\cap C_l=\emptyset$. Therefore, the number of distinct ordered
pairs that were removed from $T(\Gamma)$ can be bounded as
follows.
\begin{equation}\label{eq:pairs}
|T(\Gamma)|-|T(\Gamma ^{\shortmid})|\geq |\bigcup _{i=1}^kC_i|=\sum _{i=1}^k|C_i|= \sum _{i=1}^k 2r_i-k.
\end{equation}
By~\eqref{eq:ind}, ~\eqref{eq:edges}, ~\eqref{eq:pairs} we get
\begin{equation}
|T(\Gamma)|\geq |T(\Gamma ^{\shortmid})|+ \sum _{i=1}^k 2r_i-k\geq
|E^{\shortmid}|+(|E|-|E^{\shortmid}|)=|E|.
\end{equation}
\textbf{Case 2:} Suppose that there exist $1\leq i,j\leq k$ such
that $(i,j)\in T(\Gamma)$. We may assume that the path from $v_i$
to $v_j$ is of minimal length such that $(i,j)\in T(\Gamma)$, and
relabel the indices such that $i=1$. Let $v_m$ be a vertex such
that
$$v_m \in D^-(v_1)\cap D^-(v_j)\quad \text{or}\quad  v_m \in D^+(v_1)\cap D^+(v_j).$$
By the minimality of $k$, $m>k$. Assume that $v_m \in D^-(v_1)\cap
D^-(v_j)$ (The proof for $v_m \in D^+(v_i)\cap D^+(v_j)$ is
similar). Let $\textrm{F}:=\{1,2,\ldots ,j\}\cup \{m\}$. By
removing the vertices $\{v_s\}_{s \in \textrm{F}}$ as well as
their corresponding edges we get a new graph $\Gamma
^{\shortmid}=(V^{\shortmid},E^{\shortmid})$ with $n-j-1$ vertices,
which satisfies the theorem by the induction assumption. That is,
\begin{equation}\label{eq:ind2}
|T(\Gamma ^{\shortmid})|\geq |E^{\shortmid}|.
\end{equation}
The number of edges that were removed is
\begin{equation}\label{eq:edges2}
|E|-|E^{\shortmid}|=\sum _{t=1}^j 2r_t-(j-1)+2r_m-2.
\end{equation}
We subtract $(j-1)$, since any edge in the path from $v_1$ to
$v_j$ is counted twice. Also, we subtract $2$ since the edges
$(v_m,v_1),(v_m,v_j)$ are counted twice. Next, we count the number
of ordered pairs that were removed from $T(\Gamma)$. Again, for
every $i \in \textrm{F}$ we define the sets $C_i$. For $1\leq i<
j$
$$ C_i=\{(v_i,v)|v\in D^-(v_{i+1})\}\cup \{(v,v_i)|v\in D^-(v_{i+1})\}(\subseteq T(\Gamma)),$$
for $m$ we define
$$ C_m=\{(v_m,v)|v\in D^-(v_1)\}\cup \{(v,v_m)|v\in D^-(v_1)\}(\subseteq T(\Gamma)),$$
and for $j$
$$ C_j=\{(v_j,v)|v\in D^+(v_m)\}\cup \{(v,v_j)|v\in D^+(v_m)\}(\subseteq T(\Gamma)).$$
Similarly to~\eqref{eq:C_i} we obtain
$$|C_i|= 2r_{i+1}-1,1\leq i< k,\quad |C_m|= 2r_1-1, \quad C_j=2r_m-1.$$
By the minimality property the above sets are distinct, that is
$C_{l_1}\cap C_{l_2}=\emptyset$ for any $l_1 \neq l_2 \in
\textrm{F}.$ Therefore, the number of distinct ordered pairs that
were removed from $T(\Gamma)$ can be bounded as follows.
\begin{equation}\label{eq:pairs2}
|T(\Gamma)|-|T(\Gamma ^{\shortmid})|\geq |\bigcup _{i\in F}C_i|=\sum _{i\in F}|C_i|=\sum _{t=1}^j 2r_t+2r_m-(j+1).
\end{equation}
By~\eqref{eq:ind2}, ~\eqref{eq:edges2}, ~\eqref{eq:pairs2} we get
\begin{equation}
|T(\Gamma)|\geq |T(\Gamma ^{\shortmid})|+ \sum _{t=i}^j 2r_t+2r_m -(j+1)\geq
|E^{\shortmid}|+(|E|-|E^{\shortmid}|)=|E|.
\end{equation}
\underline{\textbf{$\Gamma$ is non-Eulerian}.}\\
First, we show that in this case there exists $(v_j,v_i)\in E$ such that
$$|D^-(v_i)|> |D^+(v_i)|\quad \text{and} \quad |D^-(v_j)|\leq |D^+(v_j)|.$$
We write $V$ as a disjoint union, $V=V_1\cup V_2$ where
$$V_1=\{v\in V||D^-(v)|> |D^+(v)|\}$$
$$V_2=\{v\in V||D^-(v)|\leq |D^+(v)|\}.$$
By the non-Eulerian property $V_1$ is not empty. Obviously, there
must be $v_j\in V_2 ,v_i\in V_1$ such that $(v_j,v_i)\in E$. Let
$v_i,v_j$ be as above and set
$$|D^-(v_i)|=i_1,|D^+(v_i)|=i_2,|D^-(v_j)|=j_1,|D^+(v_j)|=j_2.$$
Then
\begin{equation}\label{eq:notE}
i_1>i_2 \quad ,\quad j_1\leq j_2.
\end{equation}
By removing the vertices $v_i,v_j$ as well as their corresponding
edges we get a new graph $\Gamma
^{\shortmid}=(V^{\shortmid},E^{\shortmid})$ with $n-2$ vertices,
which satisfies the theorem by the induction assumption. That is,
\begin{equation}\label{eq:ind3}
|T(\Gamma ^{\shortmid})|\geq |E^{\shortmid}|.
\end{equation}
Again, let
$$ C_j=\{(v_j,v)|v\in D^-(v_i)\}\cup \{(v,v_j)|v\in D^-(v_i)\}(\subseteq T(\Gamma)),$$
and
$$ C_i=\{(v_i,v)|v\in D^+(v_j)\}\cup \{(v,v_i)|v\in D^+(v_j)\}(\subseteq T(\Gamma)).$$
The cardinality of $C_i$ is $2i_1-1$, and the cardinality of $C_j$
is $2j_2-1$. Since there are no loops in $\Gamma$, then clearly
$C_i\cap C_j=\emptyset$. Therefore, the number of distinct ordered
pairs that were removed from $T(\Gamma)$ can be bounded as
follows,
\begin{equation}\label{eq:ind4}
|T(\Gamma)|-|T(\Gamma^{\shortmid})|\geq |C_i\cup C_j|=|C_i|+|C_j|=2i_1+2j_2-2.
\end{equation}
On the other hand, by~\eqref{eq:notE} the number of edges that
were removed is bounded as follows,
\begin{equation}\label{eq:ind5}
|E|-|E^{\shortmid}|= i_1+i_2+j_1+j_2-1\leq 2i_1+2j_2-2.
\end{equation}
By~\eqref{eq:ind3},~\eqref{eq:ind4} and~\eqref{eq:ind5} we get
\begin{flalign*}
&&|T(\Gamma)|\geq |T({\Gamma^{\shortmid}})|+2i_1+2j_2-2\geq
|E^{\shortmid}|+(|E|-|E^{\shortmid}|)=|E|.&& \qed
\end{flalign*}
For an undirected graph $\Gamma$, let
$$\gamma _k=\{(i,j)|\text{there exists a path of length } k \text{ between } v_i \text{ and } v_j \}.$$
By interpreting an undirected graph as a digraph, Theorem A yields the
following corollary.
\begin{corollary}
Let $\Gamma$ be a finite undirected graph without parallel edges.
Then, with the above notations,
$$|\gamma _2|\geq |\gamma _1|.$$
\end{corollary}
Another consequence of Theorem A is as follows.
Denote the non-negative real numbers by $\mathbb{R}^+$.
For $A\in M_n(\mathbb{R}^+)$, let
$$\text{Supp}(A)=\{(i,j)\in A|a_{ij}\neq 0\}.$$
Any $A\in M_n(\mathbb{R}^+)$ determines a digraph $\Gamma=(V,E)$ in the following way:
$$V=\{v_1,v_2,\ldots ,v_n\}$$
$$(v_i,v_j)\in E \Leftrightarrow a_{ij}\neq 0.$$
Clearly $|\text{supp}(A)|=|E|$. Notice that
$(AA^t+A^tA)_{ij}\neq 0 $ if and only if the vertices $v_i,v_j$ are mutually neighbored.
\begin{corollary}
For any $A\in M_n(\mathbb{R}^+)$
$$|\text{Supp}(AA^t+A^tA)|\geq |\text{Supp}(A)|.$$
\end{corollary}
\section{gradings}\label{gradings}
In order to deduce Theorem B from Theorem A we need the following classification theorem.
\begin{theorem}\cite[see Theorem 3]{MR2488221}\label{th:zaicev2}
Let~\eqref{eq:algebragrading} be a $G$-simple grading of a complex algebra $\Lambda$. Then there exists
an $n$-tuple $(g_1,g_2,\ldots ,g_n)\in G^n$ and a subgroup $H\leq G$, where
dim$_{\mathbb{C}}(\Lambda)=n^2\cdot |H|$ such that for any $g\in G$
\begin{equation}\label{eq:induce}
\text{dim}_{\mathbb{C}}(\Lambda_g)=|\{(g_i,h,g_j)|h\in H, g_i^{-1}hg_j=g\}|.
\end{equation}
\end{theorem}
Simply-graded algebras are vastly investigated, e.g. see
\cite{giambruno, aljadeff2011simple, MR2226177, MR1941224,
MR2488221, das1999, MR2928456, MR1863551}.
Let~\eqref{eq:algebragrading} be a $G$-simple grading. Let
$(g_1,g_2,\ldots ,g_n)\in G^n$ and $H\leq G$ be the corresponding
$n$-tuple and the corresponding subgroup provided by
Theorem~\ref{th:zaicev2}. For any $g\in G$ we associate a digraph
$\Gamma _g=(V_g,E_g)$ in the following way.
$$V_g=\{v_1,v_2,\ldots ,v_n\},$$
$$E_g=\{(v_i,v_j)|\exists h\in H,\quad g_i^{-1}hg_j=g\}.$$
Notice that for any $1\leq i,j\leq n$, if $$g_i^{-1}h_1g_j=g_i^{-1}h_2g_j$$ then $h_1=h_2$.
Therefore, by~\eqref{eq:induce}
\begin{equation}\label{eq:grad}
\text{dim}_{\mathbb{C}}(\Lambda_g)=|E_g|.
\end{equation}
\begin{remark}
\cite[Theorem 3]{MR2488221} describes a way to decompose any
simple $G$-graded complex algebra to \textit{fine} and
\textit{elementary} gradings. This decomposition is not unique.
However, when given a $G$-simple
grading~\eqref{eq:algebragrading}, one can show by
\cite[Proposition 3.1]{aljadeff2011simple} and by using the
``moves" described in \cite[Lemma 1.3]{aljadeff2011simple} that
for any $g\in G$ the associated digraph $\Gamma _g$ is determined
up to a graph isomorphism.
\end{remark}
\begin{center}
Proof of Theorem B.
\end{center}
By~\eqref{eq:grad} we need to show that $|E_e|\geq |E_g|$ for any $g\in G$.
Applying Theorem A on the graph $\Gamma_g$ we obtain
\begin{equation}\label{eq:grad2}
|T({\Gamma _g})|\geq |E_g|.
\end{equation}
Now, we show that if a pair $(i,j)\in T(\Gamma _g)$  then $(v_i,v_j)\in E_e$.
Indeed, if there exists $v_k\in D^+(v_i)\cap D^+(v_j)$ (the case where $v_k\in D^-(v_i)\cap D^-(v_j)$ is similar),
then there exist $h_1,h_2$ such that
\begin{equation}\label{eq:distinct}
g_i^{-1}h_1g_k=g=g_j^{-1}h_2g_k.
\end{equation}
By~\eqref{eq:distinct}
$$(g_i^{-1}h_1g_k)(g_j^{-1}h_2g_k)^{-1}=e$$
and hence $(v_i,v_j)\in E_e.$
As a consequence we get that
\begin{equation}\label{eq:grad1}
|E_e|\geq |T(\Gamma _g)|.
\end{equation}
Therefore, by~\eqref{eq:grad2} and~\eqref{eq:grad1}
we get that for any $g\in G$
\begin{flalign*}
&& |E_e|\geq |E_g|. && \qed
\end{flalign*}

\begin{remark}
The dimension of the trivial component is not necessary maximal
when the algebra is not simply-graded. For example, consider the
natural $\mathbb{Z}$-grading of a polynomial ring with more the
one indeterminate. In this case, the trivial component is one
dimensional, whereas the other components have strictly larger
dimensions.
\end{remark}

\addcontentsline{toc}{chapter}{Bibliography}


\begin{thebibliography}{1}

\bibitem{giambruno}
A.~Aljadeff, E.~Giambruno.
\newblock Multialternating graded polynomials and growth of polynomial
  identities.
\newblock {\em Proc. Amer. Math. Soc. (in press)}.

\bibitem{aljadeff2011simple}
E.~Aljadeff and D.~Haile.
\newblock Simple g-graded algebras and their polynomial identities.
\newblock {\em Trans. Amer. Math. Soc. (in press)}.

\bibitem{MR2226177}
E.~Aljadeff, D.~Haile, and M.~Natapov.
\newblock On fine gradings on central simple algebras.
\newblock In {\em Groups, rings and group rings}, volume 248 of {\em Lect.
  Notes Pure Appl. Math.}, pages 1--9. Chapman \& Hall/CRC, Boca Raton, FL,
  2006.

\bibitem{MR1941224}
Y.~A. Bahturin and M.~V. Zaicev.
\newblock Group gradings on matrix algebras.
\newblock {\em Canad. Math. Bull.}, 45(4):499--508, 2002.

\bibitem{MR2488221}
Y.~A. Bahturin, M.~V. Zaicev, and S.~K. Sehgal.
\newblock Finite-dimensional simple graded algebras.
\newblock {\em Mat. Sb.}, 199(7):21--40, 2008.

\bibitem{das1999}
S.~D{\u{a}}sc{\u{a}}lescu, B.~Ion, C.~N{\u{a}}st{\u{a}}sescu, and
  J.~Rios~Montes.
\newblock Group gradings on full matrix rings.
\newblock {\em J. Algebra}, 220(2):709--728, 1999.

\bibitem{MR2928456}
D.~Haile and M.~Natapov.
\newblock A graph theoretic approach to graded identities for matrices.
\newblock {\em J. Algebra}, 365:147--162, 2012.

\bibitem{MR1863551}
S.~K. Sehgal and M.~V. Zaicev.
\newblock Finite gradings of simple {A}rtinian rings.
\newblock {\em Vestnik Moskov. Univ. Ser. I Mat. Mekh.}, (3):21--24, 77, 2001.

\end{thebibliography}

\end{document}